\documentclass[a4paper,12pt]{article}
\usepackage{amsmath,amssymb,amsthm}
\usepackage{url}
\usepackage{graphicx}
\newtheorem{thm}{Theorem}[section]

\newtheorem{example}[thm]{Example}
\newtheorem{definition}[thm]{Definition}

\numberwithin{equation}{section}

\def\D{\mathcal D}

\title{Forensic identification: the Island Problem and its generalisations}

\author{Klaas Slooten\footnote{Netherlands Forensic Institute, P.O. Box 24044, 2490 AA The Hague, The Netherlands;
k.slooten@nfi.minjus.nl} \,\, and Ronald Meester\footnote{VU University Amsterdam, De Boelelaan 1081, 1081 HV Amsterdam, The
Netherlands; rmeester@few.vu.nl}}

\begin{document}

\maketitle

\begin{abstract}
In forensics it is a classical problem to determine, when a suspect $S$ shares a property $\Gamma$ with a criminal $C$, the probability that $S=C$. In this paper we give a detailed account of this problem in various degrees of generality. We start with the classical case where the probability of having $\Gamma$, as well as the a priori probability of being the criminal, is the same for all individuals. We then generalize the solution to deal with heterogeneous populations, biased search procedures for the suspect, $\Gamma$-correlations, uncertainty about the subpopulation of the criminal and the suspect, and uncertainty about the $\Gamma$-frequencies. We also consider the effect of the way the search for $S$ is conducted, in particular when this is done by a database search. A returning theme is that we show that conditioning is of importance when one wants to quantify the ``weight'' of the evidence by a likelihood ratio. Apart from these mathematical issues, we also discuss the practical problems in applying these issues to the legal process. The posterior probabilities of $C=S$ are typically the same for all reasonable choices of the hypotheses, but this is not the whole story. The legal process might force one to dismiss certain hypotheses, for instance when the relevant likelihood ratio depends on prior probabilities. We discuss this and related issues as well. As such, the paper is relevant both from a theoretical and from an applied point of view.   
\end{abstract}


\section{Introduction}
In 1968, a couple stood to trial in a notorious case, known as ``People of the State of California vs.\ Collins". The pair had been arrested since it matched eye-witness descriptions. It was estimated by the prosecution that only one in twelve million couples would match this description. The jury were invited to consider the probability that the accused pair were innocent and returned a verdict of guilty.

Later, the verdict was overthrown, essentially because of the flaws in the statistical reasoning. The case sparked interest in the abstraction of this problem, which became known as the island problem, following terminology introduced by Eggleston \cite{Eggleston}.
Its formulation is the following. A crime has been committed by an unknown member of a population of $N+1$ individuals. It is known that the criminal has 
a certain property $\Gamma$. Each individual has $\Gamma$ (independently) with probability $p$. A random member of the population is tested and observed to have $\Gamma$. What is the probability that it is the criminal? 

This problem has been quite extensively studied in the literature. For example, Balding and Donnelly \cite{BD} give a detailed account of the island problem as well as of its generalization to inhomogeneous populations or (alternatively) uncertainty about $p$. They also discuss the effects of a database search or a sequential search (i.e., a search which stops when the first $\Gamma$-bearer is found). Dawid and Mortera have studied the generalization of the island problem to the case where the evidence may be unreliable \cite{DM, DM2}.

The current paper is expository in the sense that some of the above mentioned results are reproduced - albeit
presented in a somewhat different way - and a research article in the sense that we consider generalizations which to our knowledge have not appeared elsewhere. Apart from the expository versus research nature, there is another duality in this paper, namely the distinction between the purely mathematical view versus a more applied viewpoint, and we elaborate on this issue first.

Most texts focus on the ``likelihood ratio", the quantity that transforms ``prior" odds of guilt, that is, before seeing the evidence, into ``posterior" odds after seeing the evidence. There is good reason to do so. Indeed, the likelihood ratio is often viewed as the weight of the evidence - it is therefore the quantity of interest for a forensic lab, which is unable or not allowed to compute prior (or posterior, for that matter) odds, this being the domain of the court. However, this already implies a first question. Which part of the available data should be seen as the evidence, and which part is ``just" background information? In other words: which evidence do we consider and what is the context? Indeed,  
the weight of the evidence, that is, the value of the likelihood ratio, sometimes depends on which of the available information is regarded as background information or as evidence (and of course also on the propositions that one is interested in proving). 
From a purely mathematical point of view, concentrating on the``posterior probabilities", that is, the probability that a suspect is guilty, given background information and/or evidence, settles the issue. Indeed, it is well known (\cite{MS2}) that the posterior probabilities are invariant under different choices of the hypotheses as long as they are ``conditionally equivalent given the data". Hence, from a purely mathematical point of view, the situation is quite clear, and one should concentrate on the posterior probabilities rather than on the likelihood ratios. 

However, from a legal perspective things are not so simple. The likelihood ratio is, as mentioned earlier, supposed to be in the domain of the statistical expert, but what if this likelihood ratio involves prior probabilities itself? We will see concrete examples of this in this article, and in these cases the classical point of view (likelihood ratio is for the expert, the rest is for the court) does not seem to immediately apply. If we have the choice among various likelihood ratios, are there reasons to prefer one over the other? Also this question will be addressed in particular cases in this paper.       

For the island problem, the above discussion is relevant as soon as the population has subpopulations, each with their own $\Gamma$-frequency. In that case, considering the information that the criminal has $\Gamma$ as information on the one hand or as evidence on the other, leads to different likelihood ratios, but the posterior odds are (of course) the same. 
We will go into this phenomenon in detail, considering subpopulations simultaneously with uncertainty about to which subpopulation the criminal and the suspect belong, together with uncertainty about the $\Gamma$-frequencies in each of the subpopulations. Another possibility which we will consider is that of $\Gamma$-correlation or a biased search (i.e., the choice of suspect depends on the true identity of the criminal). 

The outline of this paper is as follows. In Section 2, we review the classical island problem. We then consider in Section 3 the effect of having a biased search protocol, and of having $\Gamma$-correlations; we show that these two different types of having dependencies are strongly related to each other. In Section 4, we treat the case where the population is a disjoint union of subpopulations, each with their own $\Gamma$-frequency and prior probability of having issued the criminal. In Section 5, we consider the effect of uncertainty of the $\Gamma$-frequencies, both in a homogeneous and heterogeneous population. In addition, we investigate the effect on the likelihood ratio of uncertainty about the criminal's and the suspect's subpopulations. Section 6 deals with the case in which a suspect is found through a match in 
a database. Finally, in Section 7 we present a significant number of numerical examples.

We have tried to include all details of the computations, but at the same time to state our conclusions in a non-technical
and accessible way. Our main conclusions can be recognized in the text as bulleted ($\bullet$) lists. As such, we hope that our
contribution is interesting and useful both for mathematicians, forensic scientists and legal representatives.

\section{The classical case} 
\label{classical}
Our starting point is a collection $X$ of $N+1$ individuals. 
All forthcoming random variables
are defined on a (non-specified) probability space with probability measure $P$. The random variables $C$ and $S$
take values in $X$ and represent the criminal and the suspect respectively. Furthermore, we have a characteristic $\Gamma$, for which we introduce indicator random variables $\Gamma_x$, taking value 1 if $x \in X$ has the characteristic $\Gamma$ and 0 otherwise. The $\Gamma_x$ are independent of $(S,C)$ in the sense that $P(\Gamma_x=1\mid C=y,S=z)=P(\Gamma_x=1)$ for all $x,y,z$.
The number of $\Gamma$-bearers is written as $U=\sum_{x \in X}\Gamma_x$.

We are primarily interested in the conditional probability
$$
P(C=s\mid S=s,\Gamma_C=\Gamma_s=1);
$$ 
often we follow the habit of stating the so-called {\em posterior odds} in favour of guilt, that is,
\begin{equation}\label{target}
\frac{P(C=s\mid S=s,\Gamma_C=\Gamma_s=1)}{P(C\neq s\mid S=s,\Gamma_C=\Gamma_s=1)}.
\end{equation}
Since we will often be working conditional on $\{S=s\}$ we introduce the notation $$P_s(\cdot)=P(\cdot\mid S=s).$$
We define the events $I:=\{\Gamma_C=1\}$, $G:=\{S=C\}$, $E:=\{\Gamma_S=1\}$, $E_x:=\{\Gamma_x=1\}$ and $G_x=\{x=C\}$. We will sometimes refer to the event $I$ (or similar events) as ``information", and to $E$ (or similar events) as ``evidence"; this is just colloquial use of language, and sometimes we will view $I$ as part of the evidence. 

\noindent With this notation, \eqref{target} reads
\[ \frac{P_s(G\mid I,E)}{P_s(G^c\mid I,E)},\]
which can be rewritten in two different ways, namely 
\begin{equation}
\label{target2} 
\frac{P_s(G\mid I,E)}{P_s(G^c\mid I,E)} =
\frac{P_s(E\mid G,I)}{P_s(E\mid G^c,I)} \cdot \frac{P_s(G\mid I)}{P_s(G^c\mid I)}
\end{equation}
or
\begin{equation}
\label{jawel1}
\frac{P_s(G\mid I,E)}{P_s(G^c\mid I,E)} =\frac{P_s(I, E\mid G)}{P_s(I, E\mid G^c)} \cdot \frac{P_s(G)}{P_s(G^c)}.
\end{equation}
The left hand side of these equations is called the {\em posterior odds}. In (\ref{target2}), we arrive at the posterior
odds by ``starting out" with background information $I$ via the quotient $P_s(G \mid I)/P_s(G^c\mid I),$ called the {\em prior odds}. These prior odds are transformed into the posterior odds by multiplication with $P_s(E\mid G, I)/P_s(E\mid G^c, I)$. This latter quotient is called the {\em likelihood ratio} and is supposed to be a measure of the strength of the evidence $E$. 
On the other hand, in (\ref{jawel1}) we ``start out" from prior odds $P_s(G)/P_s(G^c)$, that is, we interpreted both $I$ and $E$ as evidence. The likelihood ratio in that case is $P_s(I, E\mid G)/P_s(I, E\mid G^c)$ and measures the ``combined" strength of the evidence $I$ and $E$.

In this section treating the classical case, we assume that $C$ and $S$ are independent and that $C$ is uniformly
distributed on $X$. Furthermore, the $\Gamma_x$ are independent and identically Bernoulli distributed with 
success probability $p$. These assumptions are not without problems when applied to concrete legal cases. The assumption that $C$ is uniformly distributed means that we a priori regard each member of the population equally likely to be the criminal. It is probably the case that computations based on this assumption cannot be used as legal evidence. However, many of the computations below can also be done with other choices for the distribution of $C$. Having a particular choice in mind does allow us to compare various formulas in a meaningful way. 
The independence and equidistribution of the $\Gamma_x$ will be relaxed later on in this paper, in various ways: one can consider subpopulations with different frequencies, allow dependencies between the $\Gamma_x$ or incorporate uncertainty in the 
probability $p$. Also the independence between $C$ and $S$ will be relaxed later on. 

The outcomes in the current section do not depend on the particular $s$ we condition on, but for the sake of consistency, we do write $P_s$ instead of $P$. We abbreviate $E:=E_S$. The independence between $S$ and $C$ now implies that $P_s(G)=1/(N+1)$.
Both likelihood ratios in (\ref{target2}) and (\ref{jawel1}) are equal to $1/p$.
It easily follows that
\begin{equation}\label{sol}
\frac{P_s(G \mid I, E)}{P_s(G^c \mid I, E)}=\frac{1}{p}\cdot\frac{P_s(G\mid I)}{P_s(G^c\mid I)}=\frac{1}{p}\cdot\frac{P_s(G)}{P_s(G^c)}=\frac{1}{Np}.
\end{equation} 
In this case it does not really matter which viewpoint one takes: the likelihood is a function of $p$ alone, and does not involve any prior knowledge. Of course, as mentioned before, in a legal setting it is not clear that uniform priors are acceptable or useful, and starting from other prior probabilities is of course possible in this framework.

In the next two subsections we will examine (for the classical case) how $P_s(G\mid I,E)$ is related to the random 
variable $U$. It turns out that we may express $P_s(G\mid I, E)$ both as the inverse of the expectation of $U$ and as the expectation of the inverse of $U$, as long as we condition correctly.

\subsection{Expected number of $\Gamma$-bearers.} Before anyone is tested for $\Gamma$, $U$ has a ${\rm Bin}(N+1,p)$-distribution. When the crime is committed and it is observed that the criminal has $\Gamma$, we condition on $\Gamma_C=1$ and obtain
\begin{eqnarray*} P_s(U=k+1 \mid I)&=&\frac{P_s(I \mid U=k+1)P_s(U=k+1)}{P_s(I)} \\ &=& \frac{\frac{k+1}{N+1}{N+1 \choose k+1}p^{k+1}(1-p)^{N-k}}{p}\\ &=& {N \choose k}p^k(1-p)^{N-k}.\end{eqnarray*}
It follows that the probability that $U=k+1$, given $I$, is equal to the probability that a random variable with a ${\rm Bin}(N,p)$-distribution takes value $k$, i.e., $U\mid I$ is distributed as $1+{\rm Bin}(N,p)$. Hence, writing ${\rm E}_s$ for expectation with respect to $P_s$, we have
\[ {\rm E}_s(U\mid I)=1+Np.\]
Thus, the posterior probability of guilt is given by the inverse of the expected number of $\Gamma$-bearers, where this expectation takes into account that there is a specific individual - the criminal - who has $\Gamma$:
\begin{equation}\label{opleil} P_s(G \mid E, I)=\frac{1}{{\rm E}_s(U\mid I)}.\end{equation}
Intuitively this makes sense: the criminal is a $\Gamma$-bearer, any one of the $\Gamma$-bearers is equally likely to be the criminal, and we have found one of them. So we have to compute the expected number of $\Gamma$-bearers, given the knowledge that $C$ is one of them.

\subsection{Expected inverse number of $\Gamma$-bearers}\label{effects} As we have seen, $U\mid I$ is distributed as $1+{\rm Bin}(N,p)$. Therefore, one expects ${\rm E}_s(U\mid I)=1+Np$ bearers of $\Gamma$. If we in addition also condition on $\Gamma_S=1$, we compute 
\begin{eqnarray*} P_s(U=k \mid E, I)&=& \frac{P_s(E\mid U=k, I)P_s(U=k\mid I)}{P_s(E \mid I)}\\&=&\frac{\frac{k}{N+1}P_s(U=k\mid I)}{\frac{1+Np}{N+1}}\\&=&\frac{k}{1+Np}P_s(U=k \mid I).\end{eqnarray*}
We use this calculation to obtain: 
\begin{eqnarray} {\rm E}_s(U^{-1}\mid E, I) &=& \sum_{k=1}^{N+1} \frac{1}{k}P_s(U=k \mid E, I) \\ &=& \sum_{k=1}^{N+1}\frac{1}{k}\frac{k}{1+Np}P_s(U=k\mid I)\\&=&\label{opleil2} \frac{1}{1+Np}.\end{eqnarray}
Summarizing,
\begin{equation} \label{solution} P_s(G\mid E,I)=({\rm E}_s(U\mid I))^{-1}={\rm E}_s(U^{-1}\mid I,E).\end{equation}
So $P_s(G\mid I, E)$ is in fact also equal to the expectation of $U^{-1}$, however, not of $U\mid I$ but of $U\mid I,E$. 
This can be understood in an intuitive way: both $S$ and $C$ have $\Gamma$, they have been sampled with replacement, so the probability that they are equal is the inverse of the number of $\Gamma$-bearers. This number is unknown, so we have to take expectations, given knowledge of $S$ and $C$.

When we compare this explanation with the one of \eqref{opleil}, we see the importance of careful conditioning.

\subsection{Effect of a search, Yellin's formula} So far, $S$ and $\Gamma_S$ were supposed to be independent of
each other. In this subsection, we consider a different situation. The random variable $C$ representing the criminal is still supposed to be
uniformly distributed, but the definition of $S$ is different: we repeatedly select from $X$ - with or without
replacement - until a $\Gamma$-bearer is found, without keeping any records on the search itself, such as its duration. The $\Gamma$-bearer found
this way is denoted by $S$; if there is no $\Gamma$-bearer in the population, we set $S=*$, and define $\Gamma_*=0$. As before we write  $E= \{\Gamma_S=1\}$ 
and note that in this situation $I \subseteq E$. 

As above, we are interested in $P_s(G\mid E,I)$ which, since $I \subseteq E$, reduces to $P_s(G\mid I)$, and this conditional probability is easy to compute:
\begin{eqnarray}
\nonumber P_s(G\mid E,I) &=&  P_s(G\mid I) = \sum_{k=0}^N k^{-1}P_s(U=k\mid I) \\\label{yellin}
&=&  E_s(U^{-1}\mid I). 
\end{eqnarray}
This formula was published by Yellin in \cite{Yellin} as the solution to this version of the island problem with a search. Sometimes, however, it is incorrectly quoted in the literature (e.g. in \cite{BD}) as an incorrect solution to the island problem without search as we have discussed it. 

\subsection{Conclusions}\label{concbasic}
\begin{itemize}
\item The classical version of the island problem is not difficult to solve, but the relation between
the probability of guilt and the expected number of $\Gamma$-bearers is rather subtle. The basic formula is
$$
P_s(G\mid I,E)=({\rm E}_s(U\mid I))^{-1}={\rm E}_s(U^{-1}\mid I,E) =\frac{1}{1 + Np}.
$$
\item In the case of a search we have $I \subseteq E$ and this leads to
$$
P_s(G\mid E,I) =  E_s(U^{-1}\mid I).
$$
These outcomes are independent of $s$. 

\item For the value of the likelihood ratio, it does not matter whether or not one
interprets $I$ as background information or as evidence - in both cases the value is $1/p$ and this quantity does not depend on any
prior knowledge.
\item The prior odds, the likelihood ratio and (hence) the posterior odds are all independent of $s$.
\end{itemize}

\section{Dependencies}

In this section we relax the condition that the $\Gamma_x$ are independent random variables or that $S$ and $C$ are
independent. To this end, we define
 
\begin{equation}\label{gcorr} c_{x,y}=P(\Gamma_x=1\mid \Gamma_y=1), \end{equation}
\begin{equation}\label{sc}\sigma_{x,y}=P(S=x \mid C=y, I).\end{equation}

\subsection{Independent $\Gamma_x$}
First we assume that the $\Gamma_x$ are independent (not necessarily identically distributed) random variables, but $C$ and $S$ are not. This is the case, for instance, in a biased search situation. It also accounts for selection effects, where certain members of the population are more likely to become a suspect than others.  
We write $p_x$ for $P(\Gamma_x=1)$.
Now \eqref{target} becomes
\begin{eqnarray} \frac{P_s(E\mid G,I)P_s(G\mid I)}{P_s(E\mid G^c,I)P_s(G^c\mid I)} &=& \nonumber \frac{P(E\mid G,S=s,I)}{P(E\mid G^c,S=s,I)}\frac{P(G\mid S=s,I)}{P(G^c\mid S=s,I)}\\ \nonumber &=& \frac{1}{p_s}\frac{P(G, S=s\mid I)}{P(G^c, S=s\mid I)}\\ \nonumber &=& \frac{1}{p_s}\frac{P(S=s\mid C=s,I)}{P(S=s\mid C \neq s, I)}\frac{P(C=s \mid I)}{P(C \neq s\mid I)}\\ \label{oddsbiased}&=& \frac{1}{p_s} \frac{\sigma_{s,s}P(C\neq s\mid I)}{\sum_{y\neq s}\sigma_{s,y}P(C=y\mid I)}\frac{P(C=s \mid I)}{P(C \neq s \mid I)}
\end{eqnarray}

In this last expression \eqref{oddsbiased}, the first term $1/p_s$ is the likelihood ratio in case of a search such that $\sigma_{s,s}=\sigma_{s,y}$ for all $y \neq s$, i.e., such that the probability of selecting $s$ is independent of $C$. In particular, this holds for a search where $S$ is uniformly random but other distributions of $(S,C)$ may also satisfy this criterion.

The middle term in \eqref{oddsbiased} is the term that accounts for the bias of the search, i.e., it expresses the effect of the dependence between $S$ and $C$ in the case $S=s$.

The last term of \eqref{oddsbiased} is the ``prior odds'', the odds in favour of $C=s$, when $I$ is taken into account. It is of course also possible to start from ``prior odds'' $P(C=s)/P(C \neq s)$; this will yield the same posterior odds, but a different expression for the likelihood ratio. We will make this explicit for some special cases later on.

\subsection{Arbitrary $\Gamma_x$}
We now assume again that $S$ and $C$ are independent, but we drop the assumption that the $\Gamma_x$ are independent. In that case, we can write
\begin{eqnarray} \frac{P_s(G\mid I,E)}{P_s(G^c\mid I,E)}&=&\nonumber \frac{P_s(E\mid G,I)}{P_s(E\mid G^c,I)}\frac{P_s(G\mid I)}{P_s(G^c\mid I)}\\ \nonumber &=& \frac{P_s(G\mid I)}{P_s(E,G^c\mid I)}
\end{eqnarray}
Since we have assumed that the $\Gamma_i$ are independent of $S$ and $C$, we have \[ P_s(E\mid I, C=y)=P(\Gamma_s=1\mid\Gamma_y=1)=c_{s,y},\] and we continue as
\begin{eqnarray}
\frac{P_s(G\mid I)}{P_s(E,G^c\mid I)} \nonumber&=&\frac{P_s(G\mid I)}{\sum_{y \neq s}P_s(E,C=y\mid I)}\\ 
\nonumber &=& \frac{P_s(G\mid I)}{\sum_{y \neq s}P_s(E\mid C=y,I)P_s(C=y\mid I)}\\
\label{reskort}&=& \frac{P_s(G\mid I)}{\sum_{y \neq s}c_{s,y}P_s(C=y\mid I)}\\ 
\label{res}&=& \frac{1}{p_s}\frac{P_s(G^c\mid I)}{\sum_{y \neq s}\frac{c_{s,y}}{p_s}P_s(C=y\mid I)}\frac{P_s(G\mid I)}{P_s(G^c\mid I)}.
\end{eqnarray}

As for the case of a biased search, the term $1/p_s$ is the likelihood ratio that we obtain in the case where the $\Gamma$-correlations do not play a role, i.e., when $c_{s,y}=p_s$ for all $y \neq s$. The middle term, analogously to \eqref{oddsbiased},
\begin{equation}\label{corrfactor} \frac{P_s(C\neq s\mid I)}{\sum_{y \neq s}\frac{c_{s,y}}{p_s}P_s(C=y\mid I)}=\frac{P_s(C \neq s\mid I)}{\sum_{y\neq s}\frac{c_{y,s}}{p_y}P_s(C=y \mid I)} \end{equation}
accounts for the $\Gamma$-correlations, and the last term
\[ \frac{P_s(G\mid I)}{P_s(G^c\mid I)}=\frac{P(C=s\mid I)}{P(C\neq s\mid I)}\] describes the prior odds, conditional on $I=\{\Gamma_C=1\}$. If we remove this conditioning, we get
\begin{eqnarray}\label{reskort2} \frac{P_s(G\mid I,E)}{P_s(G^c\mid I,E)}&=&\frac{P_s(G)}{\sum_{y \neq s}c_{y,s}P_s(C=y)}\\\label{res2}&=&\frac{1}{p_s}\frac{P_s(C\neq s)}{\sum_{y \neq s}\frac{c_{y,s}}{p_s}P_s(C=y)}\frac{P_s(G)}{P_s(G^c)}.\end{eqnarray}

As for \eqref{res}, the last line contains three terms: the likelihood ratio $1/p_s$ in the uncorrelated case, the term due to the correlation and the prior odds. 

Finally, note that \eqref{reskort} and \eqref{reskort2} imply
\begin{equation}\label{IgeenI} \frac{P_s(G\mid I, E)}{P_s(G^c\mid I, E)}= \frac{P_s(G\mid I)}{\sum_{y \neq s}c_{s,y}P_s(C=y\mid I)}=\frac{P_s(G)}{\sum_{y \neq s}c_{y,s}P_s(C=y)}\end{equation} (or equivalently, the symmetry between the middle terms in \eqref{res} and \eqref{res2}):  the way the correlation between the $\Gamma_i$ appear in the posterior odds depends on whether or not one considers $I=\{\Gamma_C=1\}$ to be evidence, or an event upon which everything is conditional.

\subsection{Comparison of biased search and $\Gamma$-correlations}\label{equiv}When we compare the posterior odds (\ref{oddsbiased}) and (\ref{res}) of the two situations, we see that the expressions are very similar. Both have a correction factor in the denominator. In fact, when $S$ and $C$ are independent, then in \eqref{res} $P_s$ can be replaced with $P$, and the two cases reduce to each other if $\sigma_{x,y}/\sigma_{x,x}=c_{x,y}$ for all $x \neq y$. A trivial example of this is obtained when $C$ is uniform on $X$ and the $\Gamma_x$ are independent Bernoulli random variables. More generally, every case of a biased search without $\Gamma$-correlations where the correlation coefficients between criminal and suspect are such that $0 \leq p_y \frac{\sigma_{x,x}}{\sigma_{x,y}} \leq 1$ is equivalent (as far as the probability of guilt is considered) to a case where the search is unbiased but the $\Gamma_x$ are correlated with coefficients $c_{x,y}=p_y\frac{\sigma_{x,x}}{\sigma_{x,y}}$.

\section{Heterogeneous populations}\label{hetero} 
In this section we consider the situation where the population consists of several subpopulations, each with their own $\Gamma$-frequency and each with their own probability of containing the criminal. 
To model this, we write $X$ as a disjoint union of subpopulations $X_i$:
\begin{equation}
\label{subpops} 
X = X_1 \cup \dots \cup X_m, 
\end{equation}
with $X_i \cap X_j =\emptyset$ whenever $i \neq j$.
If $x \in X_i$, we say that $x$ is in subpopulation $i$ and write $i=X(x)$. Let $N_i=|X_i|$ be the size of subpopulation $X_i$.
We write $N_x=N_i$ if $i=X(x)$.
Let \begin{equation}\label{defbeta} P(C \in X_i)=\beta_i, \end{equation} 
where the $\beta_i$'s are positive and satisfy $\sum_{i=1}^m \beta_i=1$. 
We assume that the random variables $\Gamma_x$ are independent Bernoulli variables with probability of success $p_{X(x)}$; hence they
are not identically distributed as their distribution varies for different subpopulations.

\subsection{Posterior probability of guilt} It follows from the above that we have $c_{x,y}=p_x$ for all $x,y \in X$. Therefore, it follows from \eqref{res} and \eqref{reskort2} that
\begin{equation}\label{solheteropop} \frac{P_s(G\mid I,E)}{P_s(G^c\mid I,E)}=\frac{1}{p_s}\frac{P_s(G\mid I)}{P_s(G^c\mid I)}=\frac{P_s(G^c)}{\sum_{i=1}p_i\beta_i-p_sP_s(G)}\frac{P_s(G)}{P_s(G^c)}\end{equation}

We can work this out in more detail in the case where $S$ and $C$ are independent and $C$ is uniform on subpopulations:
\begin{equation}
\label{un}
      P(C=x\mid  C \in X(x))= 1/N_x.
\end{equation}
This assumption is not a restriction, since we assume that all $\Gamma_{x}$ are independent. It is always possible to split up the population into parts such that the $\Gamma_x$ are i.i.d. on the parts and \eqref{un} holds (a trivial decomposition would be into singletons). 

First, we define $\alpha_i$ to be the probability that $C \in X_i$, given that $C$ has $\Gamma$:
\begin{equation}\label{alphabeta} \alpha_i = P(C \in X_i \mid I)=\frac{P(I \mid C \in X_i)P(C \in X_i)}{P(I)} 
=\frac{p_i\beta_i}{\sum_{j=1}^m p_j\beta_j}.\end{equation}
Now, $P(C=x)=\alpha_x/N_x$ and $P(C=x\mid I)=\beta_x/N_x$, and \eqref{solheteropop} can be rewritten as
\begin{equation}\label{solheteropop2} \frac{P_s(G\mid I,E)}{P_s(G^c\mid I,E)}=\frac{1}{p_s}\frac{\alpha_s}{N_s-\alpha_s}=\frac{1}{N_s\sum_{i=1}p_i\frac{\beta_i}{\beta_s}-p_s}.\end{equation}

\subsection{Likelihood ratios} 
It follows from \eqref{solheteropop2} that, whether $S$ and $C$ are independent or not, the likelihood ratio conditioned on $I$ is given by
\begin{equation}\label{LR1}  \frac{P_s(E\mid G,I)}{P_s(E\mid G^c,I)}=\frac{1}{p_s}.\end{equation}
If we assume independence of $S$ and $C$ and that $C$ restricted to each subpopulation is uniform, then we obtain
\begin{equation}
\label{LR2} \frac{P_s(I,E\mid G)}{P_s(I,E\mid G^c)}=\frac{N_s-\beta_s}{N_s\sum_{j=1}^m p_j\beta_j-p_s\beta_s}.
\end{equation}
We note two special cases. First, when $N_s$ is large which means that the prior probability of guilt for $s$ is small), (\ref{LR2}) is approximately equal to 
\begin{equation}
\label{LR3}
\frac{1}{\sum_{j= 1}^m p_j\beta_j},
\end{equation}
in which the subpopulation to which $s$ belongs plays no special role.
A second special case arises when we take $N_s=1$, and only one other subpopulation. This is the standard practice for many forensic labs: there is a default population (the local population), and only two hypotheses are considered: either $S=C$, or $C$ is from the default population. 
In that case, the likelihood ratio \eqref{LR2} is equal to  
\begin{equation}
 \label{LRdefpop}
\frac{1}{p_{def}},
\end{equation}
where $p_{def}$ is the $\Gamma$-frequency in the default population and $\beta_{def}$, the prior probability that $C$ is from the default population, is equal to $1-\beta_s$.

\subsection{Discussion}
It seems that (at least) two likelihood ratios can be used to answer the informal question ``What is the weight of the evidence that the suspect has the same characteristic as the criminal?". Contrary to the classical case described in Section \ref{classical}, the weight of the evidence depends on whether or not we consider the fact that the criminal has $\Gamma$ to be evidence or background information. 
Depending on that choice and on the prior odds on guilt for $S$, we may arrive at the reciprocal of either $p_s$, $p_{def}$, or$\sum p_j \beta_j$. These quantities may be very different. This articulates the fact that one should be very careful with the use of such
likelihood ratios, and that one should primarily be interested in posterior odds rather than in likelihood ratios. A similar warning in a different situation can be found in \cite{MS1} and \cite{MS2}. 

On the other hand, if one wants to divide the ingredients in the computation of the posterior odds into parts that are for the court to decide, and parts that are for an expert witness to provide, one faces difficulties. We will now go into these in some detail. 

\subsubsection{Choice of evidence} The difference between the choice of conditioning on $I$ or not, is directly related to the difference between the questions ``What is the probability that $S$ has $\Gamma$, if innocent?" and ``What is the probability that $C$ has $\Gamma$, if $S$ is innocent?"; or more informally ``How else can we explain that $S$ has $\Gamma$?" versus ``How else can we explain that $C$ has $\Gamma$?" Indeed, if we consider both $I$ and $E$ as evidence to be expressed by a single likelihood ratio, then we can first consider $E$, and then $I$ given $E$. But without knowledge of $I$, the probability that $S$ has $\Gamma$ is the same under $G$ as under $G^c$, so the likelihood ratio of $I$ and $E$ together is in fact the same as the likelihood ratio of $I$, given $E$.
Thus, the issue here is that we need to decide if the fact that $C$ has $\Gamma$ counts as evidence against $S$, or not. Should the fact that $C$ has a certain characteristic count as (legal) evidence against someone, 
because he belongs to a subpopulation in which the characteristic is more common? Or do we only consider the fact that $S$ has the characteristic, {\it knowing} that $C$ has it, as evidence? It seems unlikely that an answer can be given in full generality, but it is important to realize that the value of the evidence will depend on it.

\subsubsection{Role of expert} Legal systems generally wish to make a distinction between the strength of the evidence, and the strength of the case. Ideally, the expert witness informs the court about the strength of the evidence (i.e., gives a Likelihood Ratio), and the court combines this information with its prior to draw conclusions about the strength of the case.
The prior is not discussed with, or communicated to, the expert. Hence, for this to be possible, the likelihood ratio should not depend on the prior of the court.
Looking at \eqref{LR2} however, it is apparent that this likelihood ratio does depend on the prior probabilities $\beta_i$ and on the suspect's population size $N_s$. 
The value of the legal evidence, if taken to be both $I$ and $E$, thus is a function of the prior and seems as such to be generally not admissible in court. In the special case \eqref{LRdefpop}, however, it is; but in that case we only obtain useful information if the assumption that either $S=C$, or $C$ is from the default population, is justified.

The Likelihood Ratio \eqref{LR1} does not suffer from these problems: it is a function of the suspect's subpopulation only, irrespective of any prior, on $S$ or on any other person or group.
Thus, if a court has somehow arrived at a prior probability $\alpha_s=P(C \in X_s \mid I)$, it can use the expert's information $p_s$ to proceed. But it must now be made clear to the court that there is a distinction between the priors with or without $I$ taken into account, and that to compute one from the other it also needs expert information.

\subsubsection{In practice: which likelihood ratio?}
We end this discussion by pointing out some pro's and cons of the likelihood ratios \eqref{LR1} and \eqref{LR2}. Clearly, 
\eqref{LR1} only involves the suspect. This is a conceptually satisfactory property, since it allows for a clear distinction between prior probabilities and the value of the evidence, as we have pointed out above. It may also provide a safeguard against using irrelevant information as evidence. Consider, for example, the following hypothetical scenario: at a crime scene, a hair of $C$ is found. Analysis by a forensic hair expert shows that $C$ must belong to subpopulation $X_1$. Later, a suspect $S\in X_1$ is found. From the hair a mitochondrial DNA profile is generated, and $S$'s mitochondrial DNA profile matches with it.
The court wishes to be informed about the value of that match. Clearly, it only makes sense to report $p_s$, since it is at this point already known that $S$ and $C$ are from the same subpopulation. But the DNA expert may not know this, and if it is standard procedure to report a variant of \eqref{LR2}, e.g.\ \eqref{LRdefpop}, then a profile frequency for the default, or even the world's population, could be reported.
 
On the other hand, an advantage of \eqref{LR2} is that it reduces the value of the evidence if there is a plausible alternative to $S$ for $C$: if there are other groups in which $\Gamma$ is relatively frequent, and which have a positive prior probability, then \eqref{LR2} decreases whereas \eqref{LR1} does not. But as we have seen, \eqref{LR2} can only do this because it makes use of all the prior probabilities, and as such it is likely to be inadmissible as legal evidence, especially if the court leaves the choice of prior to the expert. A possible way out would be for the expert to report all the $p_j$ separately to the court.

Of course, in practice $p_s$ may be hard for the expert to determine, because he only has data about other populations, or because it is not immediately clear to which subpopulation $S$ belongs, or even what the subpopulations themselves are. In that case, it may be practical (though potentially dangerous) to use \eqref{LRdefpop} and report $p_{def}$ (together with the hypotheses!), if it is the only statistic concerning $\Gamma$ that the expert has knowledge of.

The difference in numerical value of \eqref{LR1} and \eqref{LR2} may lead to the prosecution and defence having different preferences for the use of \eqref{LR1} or (a variant of) \eqref{LR2}. For example, if $p_s$ is much smaller than the weighted mean $\sum p_j\beta_j$, the prosecution will prefer \eqref{LR1}, but the defence will point out that in the population as a whole, there are subpopulations in which $\Gamma$ is much more common, and therefore try to persuade the court that \eqref{LR2} better reflects the value of the match.
The court should realize that both points of view can be justified: the prosecutor focuses on the suspect and comes up with the likelihood that $S$ has $\Gamma$, if not guilty; the defence focuses on $C$ and points out that $S$ need not be $C$, since there are other good candidates. The court should realize that these arguments can be both valid. 

To better understand the influence of uncertainty about the $\Gamma$-frequencies in the different populations and about the suspect's and the criminal's subpopulation, we proceed with a more detailed model involving these issues
in Section \ref{uncertainty}.

\subsection{Expected number of $\Gamma$-bearers} If we choose $\beta_i=N_i/N$ as we did in the classical case, then we can again express the posterior probability of guilt as the inverse of the expected number of $\Gamma$-bearers. 
We compute $E_s(U \mid I, C \in X_s)=\sum_{i}N_ip_i+1-p_s=N_s\sum_{i}p_i\frac{N_i}{N_s}-p_s+1$, and from  \eqref{solheteropop2} it follows that 
\begin{equation}
\label{} P_s(G\mid I,E)=\frac{1}{{\rm E}_s(U\mid I, C \in X_s)},
\end{equation}
which is the analogue of \eqref{opleil}. The reader may check that similarly,
\[ P_s(G \mid I, E)={\rm E}_s(U^{-1}\mid I, E, C \in X_s).\]
This is the analogue of \eqref{solution}. 

\subsection{Without conditioning on $S=s$}
Assume that $S$ is uniformly distributed on $X$, and suppose we do not condition on $\{S=s\}$. Concentrating on the conditional probability of $G$ we obtain
\begin{equation}
\label{bla}
P(G\mid I,E) = \sum_{s \in X} P(G\mid I,E,S=s)P(S=s\mid I,E).
\end{equation}
The first term in the summation is computed above already, so we need only to compute $P(S=s\mid I,E)$. Since information about $S$
and its $\Gamma$-status does not say anything about $\Gamma_C$, we have that
\begin{eqnarray*}
P(S=s\mid I,E) &=& P(S=s\mid \Gamma_S =1)\\
&=& \frac{P(\Gamma_S=1\mid S=s)P(S=s)}{\sum_{s\in X} P(\Gamma_s =1\mid S=s)P(S=s)}\\
&=& \frac{p_s}{\sum_{s \in X} p_s}.
\end{eqnarray*}
Hence it follows that 
$$
P(G\mid E,I)=\frac{\sum_{s \in X} p_s Z_s}{\sum_{s \in X} p_s},
$$
where
$$
Z_s =P_s(G\mid I,E) = \frac{\alpha_s}{p_s(N_s-\alpha_s) + \alpha_s}.
$$
Hence the posterior probability of guilt is a weighted average of the conditioned ones, with weights $p_s$.

\subsection{Conclusions}
\begin{itemize}
\item
The probability of guilt in this situation is equal to
$$
P_s(G \mid I, E)=\frac{\alpha_s}{p_s(N_s-\alpha_s)+\alpha_s},
$$
and this answer depends on $s$ via the frequency of $\Gamma$ in the subpopulation of $s$, the distribution of
$C$ and the size of the subpopulation of $s$. The sizes of the other subpopulations do not play a role other than in the assessment of the $\beta_i$ and thereby of 
the $\alpha_i$, i.e., in the distribution of $C$.
\item For the value of the likelihood ratio, it does matter whether or not $I$ is interpreted as background information
or evidence. For the probability of guilt this distinction is - of course - irrelevant, but we have seen that there can be reasons to have preference for a particular choice. It is preferable to use a likelihood ratio which does not involve any prior knowledge. The prior should then, in theory, be estimated by the juror.
\item The probability of guilt, conditioning only on the fact that the suspect has $\Gamma$ but not on the identity (subpopulation) of the suspect, is the weighted average of the individual conditional probabilities, with weight factors $p_s$. The sizes of the subpopulations and the distribution of $C$ do not play a role in the
weights.
\end{itemize}

\section{Uncertainty about the frequency of $\Gamma$}
\label{uncertainty} 
In this section we assume that the $\Gamma$-frequency $P(\Gamma_x=1)$ is not known with certainty. Instead, we describe the frequency
with a probability distribution.  

\subsection{Classical case}
We assume that there are no subpopulations. The random variable $C$ is uniform on $X$, and $S$ and $C$ are independent. To model the uncertainty of the $\Gamma$-frequency, we assume that there is a random variable $W$, taking values in $[0,1]$ and with density $\chi$, such that conditional on $W=r$, the $\Gamma_x$ are independent Bernoulli variables with $P(\Gamma_x=1)=r$. We let $p$ denote the expectation of $W$ and $\sigma^2$ its variance.
We again condition on $S=s$ whenever we compute odds, but all results in this section are independent of $s$.   
\begin{definition}
\label{dists}
The distribution of $W$ is called the {\em prior-to-crime} distribution and the distribution of $W$ conditioned on $I$ is called the {\em prior-to-suspect}
distribution. Finally, the distribution of $W$ conditioned on both $I$ and $E$ is called the {\em post-match} distribution. The densities
of these three random variables are denoted by $\chi$, $\chi_{I}$ and $\chi_{I,E}$ respectively.
\end{definition}
Since
\[ P(I)=\int_0^1 P(I\mid {W}=t)\chi(t)dt=\int_0^1 t\chi(t)dt=p,\]
the continuous version of Bayes' theorem implies that
\begin{equation}
\label{FFI} \chi_{I}(t)=\frac{P(I\mid {W}=t)\chi(t)}{P(I)}=\frac{t}{p}\chi(t).
\end{equation}
Furthermore, we have
\begin{equation}
\label{FIFIE2}
\chi_{I,E}(t)=\frac{1+Nt}{1+N(p+\sigma^2/p)}\chi_{I}(t).
\end{equation}
To see this, note that
\begin{equation}
\label{FIFIE}
\chi_{I,E}(t)=\frac{P(E\mid {W}=t, I)\chi_{I}(t)}{P(E\mid I)}
\end{equation}
and compute the denominator:
\begin{eqnarray}
P(E\mid I) &=& \int_0^1 P(E\mid {W}=t,I)\chi_{I}(t)dt\\
&=&\int_0^1\frac{1+Nt}{1+N}\frac{t}{p}\chi(t)dt \\ 
&=&  \frac{1}{(1+N)p}{(p+N(p^2+\sigma^2))}\\\label{EI} &=& \frac{1+N(p+\frac{\sigma^2}{p})}{1+N}.\label{expfi}
\end{eqnarray}
From this, the claim readily follows.

The expectation of $W$ given $I$ is expressed in terms of $\chi$ by
\begin{equation} {p':=\rm E}({W}\mid I)=\int_0^1 t\chi_{I}(t)dt=\int_0^1\frac{1}{p}t^2\chi(t)dt=\frac{1}{p}(p^2+\sigma^2).\end{equation}
The expected number of $\Gamma$-bearers, given $I$ is now given by
\begin{equation}
\label{ceciel}
{\rm E}({U}\mid I) = \int_0^1 {\rm E}({U} \mid I, {W}=t)\chi_{I}(t)dt = 
\int_0^1 (1+Nt)\chi_{I}(t)dt = 1+Np'.
\end{equation}
As in the classical case where $\sigma^2=0$ (cf.\ \eqref{opleil}), the inverse of this expression is equal to the posterior probability of guilt, since 
\begin{eqnarray} P_s(G \mid I, E)&=&
\int_0^1 P_s(G \mid I, E, {W}=t)\chi_{I,E}(t)dt \\ &=& \int_0^1\frac{1}{1+Nt}\frac{1+Nt}{1+N(p+\sigma^2/p)}\chi_{I}dt\\
\label{soluncert}&=& \frac{1}{1+N(p+\sigma^2/p)}=\frac{1}{1 + Np'}.
\end{eqnarray}
Since the prior probability of guilt is just $1/(N+1)$ as before, the likelihood ratio is $1/p'$. Since this likelihood ratio is not controversial in this case, we concentrate on the posterior 
probability of guilt in terms of the various conditional distributions.

As in the classical case (cf.\ \eqref{opleil2}), we also have $P_s(G\mid I,E)={\rm E}_s({U}^{-1}\mid I,E)$. Indeed, 
\begin{eqnarray}
{\rm E}_s({U}^{-1}\mid I,E)&=& \int_0^1 {\rm E}_s({U}^{-1} \mid I,E, {W}=t)\chi_{I,E}(t)dt \\
&=& \int_0^1 \frac{1}{1+Nt}\frac{1+Nt}{1+N(p+
\sigma^2)}\chi_{I}(t)dt  \\
&=&\frac{1}{1+N(p+\sigma^2/p)}. \label{derde}
\end{eqnarray}

The expectation $p'$ only depends on $\chi$ and not on the population size. This is to be expected, since learning that a (randomly
chosen) population member has $\Gamma$ is not informative about the population size. This changes when we learn $E$, the fact that a randomly selected islander has $\Gamma$ as well. Indeed, in a small population this is more likely to happen since we are more likely to accidentally select the criminal. In the extreme case where $N=0$, $E$ can not offer any new information, but for other $N$, it does. 
It follows from \eqref{FIFIE2} that
\begin{eqnarray*} 
p''&:=& {\rm E}_s({W}\mid I,E)=\int_0^1t\chi_{I,E}(t)dt\\ &=& \frac{1}{1+N(p+\sigma^2/p)}\int_0^1t(1+Nt)\frac{t}{p}\chi(t)dt \\ &=& \frac{1}{1+N(p+\sigma^2/p)}\left(p+\frac{\sigma^2}{p}+\frac{N}{p}\int_0^1t^3\chi(t)dt\right).
\end{eqnarray*}
We can also write 
\[ p''=\frac{1}{1+Np'}\left(p'+N(\sigma_{{W}\mid I}^2+p'^2)\right),\]
if we want to express $p''$ in terms of $\chi_{I}$, where $\sigma^2_{W\mid I}$ denotes the variance of $\chi_I$.
The above formula can be rewritten as
\[ p''=p'\frac{1+N(\sigma_{{W}\mid I}^2/p'+p')}{1+Np'} \geq p',\]
with equality only if $\sigma_{{W}\mid I}^2=0$ or $N=0$ (as expected, cf.\ the remark above).

It is perhaps worth mentioning that one can reconstruct $\chi_I$ from $\chi_{I,E}$ and $\chi$ from $\chi_I$. Indeed we have

\begin{equation}\label{FIEFI} 
\chi_I(t)=\frac{\chi_{I,E}(t)}{1+Nt}\left( \int_0^1 \frac{\chi_{I,E}(s)ds}{1+Ns}\right)^{-1}
\end{equation}
and
\begin{equation}
\label{tweede}
\chi(t)=\left(t\int_0^1 \frac{\chi_I(x)}{x}dx\right)^{-1}\chi_I(t).
\end{equation}  

To see this, note that from \eqref{FIFIE2} we have
\begin{equation}\label{FIFIE3} \chi_I(t)=\frac{1+Np'}{1+Nt}\chi_{I,E}(t).\end{equation}
On the other hand, it follows from  (\ref{derde}) that
\[ {\rm E}_s({U}^{-1}\mid I,E)=\int_0^1 \frac{1}{1+Ns}\chi_{I,E}(s)ds=\frac{1}{1+Np'},\]
and the first claim (\ref{FIEFI}) follows.

For (\ref{tweede}) we simply note from (\ref{FFI}) that
\begin{equation}
\label{gg}
\chi(t)=\frac{p}{t}\chi_I(t),
\end{equation}
where $p=\int_0^1 t \chi(t)dt$. Integrating this equation gives
$$
1=p\int_0^1 \frac{\chi_I(t)}{t}dt=1
$$
and this expresses $p$ in terms of $\chi_I$. Substituting this into (\ref{gg}) gives (\ref{tweede}).

As a conclusion, we have seen that
\[ p'' \geq p' \geq p,\] so one has
\[ \frac{1}{1+Np''} \leq \frac{1}{1+Np'} \leq \frac{1}{1+Np}.\]

\subsubsection{Conclusions}
\begin{itemize}
\item The basic formula of Conclusions \ref{concbasic} still holds: using (\ref{ceciel}), (\ref{soluncert}) and (\ref{derde}), we see that the probability of guilt is given by
\[ P_s(G\mid I,E)={\rm E}_s({U}^{-1}\mid I,E)={\rm E}_s({U}\mid I)^{-1}=\frac{1}{1+Np'}.\]
\item 
The conditional probability of guilt expressed in terms of $\chi$ is 
\begin{equation}\label{ptc} P_s(G\mid I, E)=\frac{1}{1+N(p+\sigma^2/p)}.\end{equation}
Therefore, ignoring the uncertainty (i.e., using $p$ instead of $p'$), is unfavourable to the suspect. If, on
the other hand, one incorrectly assumes that there is uncertainty, then this is favourable to the suspect.
\item 
The conditional probability of guilt expressed in terms of $\chi_{I}$ is
\begin{equation}\label{pte} P_s(G \mid I,E)=\frac{1}{1+Np'}.\end{equation}
In this case, the uncertainty in $\chi_{I}$ is irrelevant in the sense that its variance plays no role.
\item 
The conditional probability of guilt expressed in terms of $\chi_{I,E}$ is 
\begin{equation}\label{pcd} P_s(G \mid I,E)=\int_0^1\frac{1}{1+Nt}\chi_{I,E}(t)dt.\end{equation} 
Ignoring the uncertainty in $\chi_{I,E}$ (obtaining $P_s(G\mid I,E)=1/(1+Np'')$) would be favourable to the suspect.  
\end{itemize}

\subsection{Uncertainty about the criminal's subpopulation}\label{pops}

Suppose that, as in Section \ref{hetero}, the population is divided into subpopulations $X=X_1 \cup \dots \cup X_m$, and that $C$ has characteristic $\Gamma$. We let $W_i$ be the random variable modelling the frequency of $\Gamma$ in $X_i$. The expectation resp.\ variance of $W_i$ are denoted by $p_i$ resp.\ $\sigma^2_i$. So, if $X(x) \neq X(y)$ then $\Gamma_x$ and $\Gamma_y$ are independent, and furthermore conditional on $W_i=p_i$ the $\Gamma_x$ for $x \in X_i$ are independent Bernoulli variables with probability of success $p_i$.
We write $E=E_s$ and $G=G_s$ as before.
Contrary to the situation in \ref{subpops}, the division of $X$ into subpopulations is a real restriction: the $\Gamma_x$ are only independent between subpopulations, not within one (only exchangeable).

\subsubsection{Unconditioned on $I$} 
We first interpret $I$ as evidence, not as background information.
The posterior probability of guilt, given that $S=s  \in X_s$ and $C$ has $\Gamma$, is
\[ P_s(G\mid I,E)=P_s(G\mid C \in X_s, I,E)P_s(C \in X_s\mid I,E).\]
The first term in the right hand side equals (see (\ref{soluncert}))
\[ P_s(G \mid C \in X_s, I, E)=\frac{1}{1+(N_s-1)(p_s+\sigma^2_s/p_s)},\]
since we are now back in the setting of a homogeneous population.
The second term equals (cf. \eqref{EI},\eqref{alphabeta})
\begin{eqnarray*} P_s(C \in X_s \mid I, E)&=&\frac{P_s(E\mid C \in X_s, I)P_s(C \in X_s \mid I)}{P_s(E \mid I)},\\ &=& \frac{\frac{1+(N_s-1)(p_s+\sigma^2_s/p_s)}{N_s}\frac{p_s\beta_s}{\sum_{j=1}^mp_j\beta_j}}{P_s(E\mid I)}.\end{eqnarray*}
It remains to compute $P_s(E\mid I)$:
\begin{eqnarray*} P_s(E\mid I) &=& \sum_{j=1}^mP_s(E \mid C \in X_j)P_s(C \in X_j \mid I) \\ &=& \sum_{j=1, j \neq X(s)}^m p_s \frac{p_j\beta_j}{\sum_{k=1}^mp_k\beta_k}+\frac{1+(N_s-1)(p_s+\sigma_s^2/p_s)}{N_s}\frac{p_s\beta_s}{\sum_{k=1}^mp_k\beta_k} \\
&=& \frac{N_s\sum_{j=1,j\neq X(s)}p_s\beta_jp_j+(1+(N_s-1)(p_s+\sigma^2_s/p_s))\beta_sp_s}{N_s\sum_{k=1}^m\beta_kp_k}.\end{eqnarray*}
Putting the parts together yields
\begin{equation}\label{post} P_s(G \mid I,E)=\frac{1}{1+N_s \sum_{j=1}^m\frac{\beta_j}{\beta_s}p_j+(N_s-1)\frac{\sigma^2_s}{p_s}-p_s}.\end{equation}
This is the analogue of \eqref{solheteropop2}. 
For large $N_s$, the probability of guilt is roughly equal to
\begin{equation}
\label{post2}
P_s(G\mid I,E) \approx \frac{1}{N_s(\sum_{j=1}^m\frac{\beta_j}{\beta_s}p_j+\sigma^2_s/p_s)}.
\end{equation}
The odds on guilt are then roughly equal to
\begin{equation}\label{postodds2} \frac{P_s(G \mid I, E)}{P_s(G^c \mid I, E)}=\frac{1}{N_s(\sum_{j=1}^m\frac{\beta_j}{\beta_s}p_j+\frac{\sigma_s^2}{p_s})}=\frac{\beta_s}{N_s}\frac{1}{\sum_{j=1}^m\beta_jp_j+\beta_s\frac{\sigma_s^2}{p_s}}.\end{equation}

The weight of the evidence, the likelihood ratio, is given by
\begin{equation}
\label{LREI}
\frac{P_s(E,I\mid G)}{P_s(E,I\mid G^c)} = \frac{N_s-\beta_s}{N_s\sum_{j=1}^mp_j\beta_j+N_s(p_s'-p_s)-p_s'\beta_s},
\end{equation}
(where $p_s'=p_s+\frac{\sigma_s^2}{p_s}$ as before) which reduces to \eqref{LR2} if $p_s'=p_s$. For large populations, \eqref{LREI} is roughly equal to 
\begin{equation}\label{LR6} \frac{1}{\sum_{j=1}^m\beta_jp_j+\beta_s\sigma^2_s/p_s},
\end{equation}
as is also clear from \eqref{postodds2}. This formula is the analogue of (\ref{LR3}). The likelihood ratio \eqref{LR6} suffers from the same problem as \eqref{LR3} in the sense that the prior probabilities $\beta_s$ are needed to compute it. For the same reason as before, it is therefore highly questionable whether the expert is allowed to report \eqref{LR6} in court. Therefore, we proceed by working conditional on $I$ and see what the computations tell us there.  

\subsubsection{Conditional on $I$} 
Now, $I$ is interpreted as background information. Let, as in \eqref{alphabeta},
\[ \alpha_s=P_s(C \in X_s \mid  I)=\frac{p_s\beta_s}{\sum_{k=1}^mp_j\beta_j}.\]
Then \eqref{post} can be rewritten as
\[ P_s(G\mid I,E)=\frac{1}{\sum_{j \neq X(s)}N_sp_s\frac{\alpha_j}{\alpha_s}+1+(N_s-1)p_s'},\]
where
\[ p_i'=p_i+\frac{\sigma_i^2}{p_i}\] 
is the expectation of $W_i$ given $C \in X_i$.

Since the prior odds, conditional on $I$, in favour of guilt of $s \in X_s$ are
\[ \frac{P_s(G\mid I)}{P_s(G^c\mid I)} =\frac{\alpha_s}{N_s-\alpha_s},\]
the corresponding likelihood ratio is equal to
\begin{eqnarray*}
\frac{P_s(E\mid G,I)}{P_s(E\mid G^c,I)}& =& \frac{N_s-\alpha_s}{N_sp_s(1-\alpha_s)+\alpha_s(1+(N_s-1)p_s')} \\
&=& \frac{N_s-\alpha_s}{N_s\alpha_s(p_s'-p_s)+N_sp_s-p_s'\alpha_s}.
\end{eqnarray*}
Of course, this leads to the same posterior probability of guilt as in (\ref{post2}).
Notice that when $p_s'=p_s$, then this reduces to $1/p_s$, i.e., we retrieve \eqref{LR1}.

For large $N_s$, the likelihood ratio is roughly equal to 
\begin{equation}\label{LR4} 
\frac{1}{p_s+\alpha_s(p_s'-p_s)}=\frac{p_s}{p_s^2+\alpha_s\sigma_s^2}.
\end{equation}
This likelihood ratio also depends on prior quantities, this time on $\alpha_s$. Note however that there is a difference between \eqref{LR6} and \eqref{LR4}. The latter only depends on quantities associated to the suspect's subpopulation, whereas the former does not. In this case there is a way to deal with the problem of having a prior quantity entering the formula for the likelihood ratio. In \eqref{LR4} one can be conservative and take $\alpha_s=1$ to obtain a number which is not larger than the true likelihood ratio. In \eqref{LR6} one can of course do the same for all $\beta_j$'s but there we have the problem that we have various $\beta_j$ in the expression, and the only thing we know is that they add up to 1. Therefore, we prefer \eqref{LR4}, but the usual care must be exercised when using this likelihood ratio in court. The use of this likelihood ratio is, as always, dangerous and should involve a discussion of priors. A likelihood ratio out of context is not useful, and unfortunately, the context is rather complicated.

\subsubsection{Conclusions}
\begin{itemize} 
\item As in the case without uncertainty about the $\Gamma$-frequencies, we obtain two likelihood ratios that quantify the weight of the evidence: for large populations these are \eqref{LR6} if the evidence is taken to be $(I,E)$ and \eqref{LR4} if the evidence is taken to be only $E$. Sine \eqref{LR4} can be easily be turned into a conservative bound by setting $\alpha_s=1$, we prefer to use \eqref{LR4}, noting however that a report mentioning just the likelihood ratio without context is dangerous and potentially misleading.
\item Only the uncertainty about the frequency of $\Gamma$ in the suspect's subpopulation plays a role in the likelihood ratio and the posterior probability of guilt, the uncertainty in the other subpopulations does not. The effect of this uncertainty is weighted by the probability that the true culprit belongs to this subpopulation. 
\item As in the classical case, if one conditions on $I$ then the likelihood ratio given by \eqref{LR4} for large populations, only contains quantities
associated to the suspect's subpopulation. 
\item Contrary to the classical case, if one considers the evidence to be $(I,E)$ then in the likelihood ratio for large populations (given by \eqref{LR6}) the suspect's subpopulation plays a special role, through the uncertainty about the $\Gamma$-frequency in this population. 
\item Regardless of whether one lets the evidence be $I,E$ or only $E$, the greater the uncertainty, the lower the weight of the evidence.
\end{itemize}

\subsection{Uncertainty about the suspect's and the criminal's subpopulation}\label{uncertain}

Suppose now that it is also unknown to which subpopulation $s$ belongs. In that case we can no longer condition on $S=s$, but
we can use the results of the previous section by writing
\begin{equation}\label{splits} P(G \mid I, E)=\sum_{i=1}^m P(G \mid S \in X_i, I, E) P(S \in X_i \mid I, E).\end{equation}

We have determined the $P(G \mid S=s, I, E)$ in \eqref{post}, and it is not difficult to see
that this is equal to $P(G\mid S\in X_i, I, E)$ whenever $s \in X_i$. Hence, we only need to
compute $P(s \in X_i \mid I, E)$.

The distribution of $S$ plays a role now, and we define
\[ \epsilon_i=P(s \in X_i)\]
to be the probability that $S$ belongs to $X_i$. 
Then the a priori probability of guilt is
\[ P(G)=P(C=s)=\sum_{i=1}^mP(s \in X_i)P(C=s \mid s \in X_i)=\sum_{i=1}^m \epsilon_i\frac{\beta_i}{N_i}.\]
Recall that $\beta_i$ is the probability that $C \in X_i$ and that we assume a uniform distribution over each subpopulation.

We now compute $P(S \in X_i \mid I,E)$:
\begin{eqnarray} \nonumber P(S \in X_i\mid I,E)&=&\frac{P(E\mid S \in X_i, I)P( S\in X_i \mid I)}{P(E \mid I)}\\&=& \label{som} \frac{P(E\mid S \in X_i, I)P( S\in X_i \mid I)}{\sum_{j=1}^mP(E \mid S \in X_j, I)P(S \in X_j \mid I)}.\end{eqnarray}
It remains to compute $P(E \mid S \in X_j, I)$ and $P(S \in X_i \mid I)$.
The latter is easy: since $I$ is information about $C$ and not about $S$, we have
\[ P(S \in X_i \mid I)=\epsilon_i.\]
The former can be computed as follows:
\[ P(E \mid S \in X_i, I)=\sum_{j=1}^m P(E \mid C\in X_j, S \in X_i, I)P(C \in X_j \mid S \in X_i, I).\]
Now $P(C \in X_j \mid S \in X_i, I)$ is the probability that $C$ belongs to $X_j$, given that $S$ has been selected from $X_i$ and that $C$ has $\Gamma$. However, nothing is given about $S$'s $\Gamma$-status and therefore $S\in X_i$ can not be informative about $C$ at all, hence 
\[ P(C\in X_j \mid S \in X_i, I)=P(C\in X_j \mid I)=\frac{p_j\beta_j}{\sum_{k=1}^mp_k\beta_k}.\]
It remains to evaluate the terms $P(E \mid C \in X_j, S \in X_i, I)$. If $i \neq j$ then $S$ and $C$ belong to different populations. If $i=j$ then \eqref{EI} applies, so 
\[ P(E \mid S \in X_i, C \in X_j, I)=\begin{cases} p_i & i\neq j, \\ \frac{1+(N_i-1)(p_i+\sigma^2_i/p_i)}{N_i} & i=j.\end{cases}\]  

If we put these ingredients together, we obtain after some computations:

\[ P(E \mid S \in X_i, I)= \frac{p_i}{\sum_{k=1}^mp_k\beta_k}\left( \sum_{j=1}^m p_j\beta_j+\frac{\beta_i}{N_i}(1-p_i+(N_i-1)\sigma_i^2/p_i) \right).\]
Plugging this into \eqref{som}, we obtain
\begin{eqnarray} \nonumber P(S \in X_i \mid I, E)&=& \frac{p_i( \sum_{k=1}^m p_k\beta_k+\frac{\beta_i}{N_i}(1-p_i+(N_i-1)\sigma_i^2/p_i) )\epsilon_i}{\sum_{j=1}^m  p_j( \sum_{k=1}^m p_k\beta_k+\frac{\beta_j}{N_j}(1-p_j+(N_j-1)\sigma_j^2/p_j) )\epsilon_j} \\ \label{probsinxi} &=& \frac{p_i\epsilon_i\beta_i}{N_iP_s(G \mid I,E, S \in X_i)}\frac{1}{\sum_{j=1}^m\frac{p_j\epsilon_j\beta_j}{N_jP_s(G \mid I,E, S \in X_j)}} \end{eqnarray}
Substituting this expression into \eqref{splits}, we arrive at the posterior probability of guilt:
\begin{equation}\label{eind} P(G \mid I,E)=\frac{\sum_{i=1}^m\frac{p_i\epsilon_i\beta_i}{N_i}}{\sum_{i=1}^m\frac{p_i\epsilon_i\beta_i}{N_iP_s(G\mid I,E, S \in X_i)}}.\end{equation}

Although this is not immediately obvious in the above presentation, the expression \eqref{eind} is symmetric in $\epsilon$ and $\beta$. To show this, notice that we only have to prove it for the denominator. Denoting 
\begin{eqnarray*} f(\epsilon,\beta) &=& \sum_{i=1}^m\frac{p_i\epsilon_i\beta_i}{N_iP_s(G\mid I,E, S \in X_i)}\\ &=& \sum_{i=1}^m \frac{p_i\beta_i\epsilon_i}{N_i}(1+N_i\sum_{j=1}^mp_j\beta_j/\beta_i+(N_i-1)\sigma_i^2/p_i-p_i),
\end{eqnarray*}
we compute 
 \begin{eqnarray*} f(\epsilon,\beta)-f(\beta,\epsilon) &=& \sum_{i=1}^m p_i\beta_i\left(\sum_{j=1}^mp_j\beta_j/\beta_i-\sum_{j=1}^mp_j\epsilon_j/\epsilon_i\right) \\ &=& 
\sum_{i=1}^m p_i \left(\epsilon_i\sum_{j=1}^mp_j\beta_j-\beta_i\sum_{j=1}^mp_j\epsilon_j \right)\\ 
&=& \sum_{i,j=1}^m \left(p_i\epsilon_ip_j\beta_j-p_i\beta_ip_j\epsilon_j\right) \\
&=& 0.
\end{eqnarray*}
Intuitively, it is clear that \eqref{eind} must possess this symmetry. Indeed, we have an unknown criminal $C$ 
and a suspect $S$, both with $\Gamma$. The probability that $S=C$ depends, as far as $\epsilon$ and $\beta$ are concerned, on how they allow for $S$ and $C$ to be issued from the same subpopulation. Exchanging the distributions $\epsilon$ and $\beta$ should not make a difference.

To conclude this section we sketch the behaviour of \eqref{probsinxi} in extreme situations.

\subsubsection{Probability that $S \in X_i$ for extreme situations}
\begin{itemize}
\item
If all $\sigma^2_j=0$ and the $N_j$ are very large (compared to the $p_j^{-1}$), then \eqref{probsinxi} is approximately equal to
\[ P(S \in X_i \mid I, E) \approx \frac{p_i\epsilon_i}{\sum_{j=1}^mp_j\epsilon_j}=P(S \in X_i \mid E).\]
This is reasonable, since if $p_i$ is big compared to $1/N_i$, then it is very unlikely that $C=S$ even when $\Gamma$ is taken into account. In this case, knowing that $C$ has $\Gamma$ does not really alter our belief about $S$'s subpopulation which we have based on $E$.
\item
If all $\sigma^2_j=0$ and the $p_j$ are small compared to the $1/N_j$, then 
\begin{equation}\label{limiet} P(S \in X_i \mid I, E) \approx \frac{p_i\frac{\beta_i}{N_i}\epsilon_i}{\sum_{j=1}^m p_j\frac{\beta_j}{N_j}\epsilon_j}=\frac{1}{\sum_{j=1}^m \frac{p_j}{p_i}\frac{\beta_j}{\beta_i}\frac{N_i}{N_j}\frac{\epsilon_j}{\epsilon_i}}.\end{equation}
If $\epsilon_i=N_i/N$, then \eqref{limiet} reduces to
\begin{equation} P(S \in X_i \mid I, E) \approx \frac{p_i\beta_i}{\sum_{i=1}^m p_j\beta_j}=P(C \in X_i \mid I),\end{equation}
which is also reasonable, since for very small $\Gamma$-frequencies it is quite likely that $C=S$.

\item
If also $\beta_i=N_i/N$, then \eqref{limiet} reduces to
\begin{equation} P(S \in X_i \mid I, E) \approx \frac{N_ip_i}{\sum_{i=1}^mN_jp_j},\end{equation}
and this is also understandable: if there is no information about the identity of $C$ or $S$, then the probability that $S \in X_i$ is proportional to the expected number of $\Gamma$-bearers in that subpopulation.
\end{itemize}

\section{Database search}

In this section we suppose that there is a database $\D \subset X$ containing the $\Gamma$-status of individuals $x_1,\dots,x_n$.
After possibly renumbering, we write $X=\{x_1,\dots,x_{N+1}\}$ and let $\D=\{x_1,\dots,x_n\}$. Suppose that $\sum_{d \in D}\Gamma_d=k$, that is, there are $k$ matches in the database. Let the evidence $E_{\D}$ be given by 
$$
E_{\D}=\{\Gamma_{x_1}=\dots=\Gamma_{x_k}=1, \Gamma_{x_{k+1}}=\dots=\Gamma_{x_n}=0\}.
$$
We also assume that $P(C=x_i\mid I)=\alpha_i$
and that each individual has $\Gamma$ with probability $p$.

There are several pairs of propositions whose support by the data can be considered. These propositions all give rise to their own likelihood ratios or posterior probabilities, which has caused considerable confusion in the literature; see \cite{MS1}
for an account on this. Some of the forthcoming discussion also appears in \cite{MS1} but we recall it here for completeness. 

We will discuss three ways of looking at database matches. The most interesting case is where the database search produces a single match. Indeed, if there are no matches then the inquiry comes to an end as far as the database is concerned and if there are several matches, then it is clear that chance matches have occurred:
\begin{enumerate}
\item Database-focused: in this case, the quantity of interest is $P(C \in \D\mid E_\D, I)$, the probability that the criminal is in the database;
\item Individual-focused: in this case, the quantity of interest is $P(C=x_1 \mid E_\D, I)$, the conditional probability that $C=x_1$ supposing that $x_1$ has $\Gamma$;
\item Database effectiveness: in this case, the quantity of interest is $P(S=C\mid E_1, I)$, the probability that $S=C$ where $E_1$ denotes the event that $k=1$ (a unique match, but not specified with whom), and where $S$ is the label of the matching individual.
\end{enumerate}
\subsection{Database-focused} First, we consider the proposition, found e.g.\ in \cite{Stockmarr}, 
\[ C \in \D,\]
and its negation $C \notin \D$. The prior odds in favour of $C \in \D$ are
\[ \frac{P(C \in \D \mid I)}{P(C \notin \D \mid I)}=\frac{\alpha_1+\dots+\alpha_n}{\alpha_{n+1}+\dots+\alpha_{N+1}},\]
where $\alpha_i=P(C=x_i \mid I)$ is the probability of guilt of $x_i$, given that $C$ has $\Gamma$.
Clearly, 
\[ P(E_\D\mid C \notin \D, I)=p^k(1-p)^{n-k}.\]
Similarly, it is easy to see that
\[ P(E_\D \mid C \in \D, I)=\frac{p^{k-1}(1-p)^{n-k}(\alpha_1 + \cdots + \alpha_k)}{\alpha_1+\dots+\alpha_n},\]
and therefore the likelihood ratio of evidence $E_\D$ in favour of $C \in \D$ is equal to
\begin{equation}\label{LRDB} \frac{P(E_\D \mid C \in \D, I)}{P(E_\D \mid C \notin \D, I)}=
\frac{\alpha_1 + \cdots + \alpha_k}{p(\alpha_1+\dots+\alpha_n )}.\end{equation}
The posterior odds in favour of $C \in \D$ are 
\begin{equation}\label{CDodds} \frac{P(C \in \D \mid E_\D, I)}{P(C \notin \D \mid E_\D, I)}=
\frac{\alpha_1 + \cdots + \alpha_k}{p (\alpha_{n+1}+\dots+\alpha_{N+1})}.\end{equation}
If $k=1$, $C \in \D$ becomes logically equivalent to $C=x_1$, and we have
\begin{eqnarray}\label{k=1} 
\frac{P(C=x_1 \mid E_\D, I)}{P(C \notin \D \mid E_\D, I)}&=& \frac{P(C=x_1 \mid E_\D, I)}{P(C \neq x_1 \mid E_\D, I)}\\
&=& \frac{\alpha_1}{p(\alpha_{n+1}+\dots+\alpha_{N+1})}=\frac{1}{p}\frac{P(C=x_1 \mid I)}{P(C \notin \D \mid I)}.
\end{eqnarray}
This means that the likelihood ratio is uncontroversial and equal to $1/p$.
In fact, it is not difficult to show that \eqref{k=1} also holds when the probability of having $\Gamma$ differs among the individuals in the database. In that case, $p$ in \eqref{k=1} should be replaced with $p_1=P(\Gamma_{x_1}=1\mid I)$. Therefore, the weight of the evidence is not influenced by the presence in the database of people of different ethnic origin other than by the determination of the $\alpha_i$.

\subsection{Individual-focused} Of course, the proposition $C \in \D$ is not really of interest to a court. Rather, presented with an individual $x$ such that $\Gamma_x=1$, a court is interested in $P(C=x \mid E_\D, I)$.
Therefore, suppose as above that there are $k$ hits in the database, namely $x_1,\dots, x_k$. A computation analogous to the above one shows that the posterior odds in favour of $C=x_1$ are
\begin{equation}\label{postdat} \frac{P(C=x_1 \mid E_\D, I)}{P(C \neq x_1 \mid E_\D,I)}=
\frac{\alpha_1}{\alpha_2 + \cdots + \alpha_k + p(\alpha_{n+1} + \cdots + \alpha_{N+1})}.\end{equation}
Notice that, if $k=1$, we retrieve \eqref{k=1}, as we should.

\subsection{Database effectiveness} The most interesting case is when the database produces a unique hit. In that case, as we have seen, the posterior odds in favour of $S=C$ are given by \eqref{k=1}. In this section we investigate a related, but different probability, namely the probability that if we have a unique database hit, that it is with the true culprit. This probability represents the long term effectiveness of the database in selecting the correct individual in the cases where it produces a unique match. We let $E_1$ denote the event that there is exactly one $\Gamma$-bearer in the database, and we will calculate
\[ P(S=C \mid E_1, I),\]
where $S$ is the unique individual in the database with $\Gamma$.
To do so, we write
\[ P(S=C\mid E_1, I)=\sum_{i=1}^nP(S=C \mid \Gamma_{x_i}=1, E_1, I)P(\Gamma_{x_i}=1\mid E_1, I).\]
First notice that \eqref{k=1} gives
\[ P(S=C \mid E_1, \Gamma_{x_i}=1, I)=\frac{\alpha_i}{\alpha_i+pP(C \notin \D \mid I)},\]
and it remains to compute $P(\Gamma_{x_i}=1\mid E_1, I)$: 
\begin{eqnarray*}  P(\Gamma_{x_i}=1 \mid E_1, I)&=& P(C=x_i\mid E_1, I)+\frac{1}{n}P(C \notin \D \mid E_1, I)\\&=&\frac{P(E_1\mid C=x_i, I)P(C=x_i\mid I)+\frac{1}{n}P(E_1\mid C \notin \D \mid I)P( C \notin \D \mid I)}{P(E_1\mid I)}\\&=&\frac{P(E_1\mid C=x_i, I)P(C=x_i\mid I)+\frac{1}{n}P(E_1\mid C \notin \D, I)P( C \notin \D \mid I)}{P(E_1\mid C\in D, I)P(C \in D \mid I)+P(E_1 \mid C \notin \D, I)P(C \notin \D \mid I)}\\&=&\frac{\alpha_i+pP(C \notin \D \mid I)}{P(C \in D \mid I)+npP(C \notin \D \mid I)},  
\end{eqnarray*}
where in the last step we used that $P(E_1\mid C \in \D, I)=(1-p)^{n-1}$ and $P(E_1\mid C \notin \D,I)=p(1-p)^{n-1}$.
It follows that
\begin{eqnarray*} P(S=C\mid E_1, I)&=&\sum_{i=1}^n \frac{\alpha_i}{\alpha_i+pP(C \notin \D \mid I)}\frac{\alpha_i+pP(C \notin \D \mid I)}{P(C \in D\mid I)+npP(C \notin \D \mid I)}\\ &=&\frac{P(C \in \D \mid I)}{P(C \in D \mid I)+npP(C \notin \D \mid I)}.\end{eqnarray*}
which can also be written in odds form:
\begin{equation}\label{effodds} \frac{P(C \in \D \mid E_1, I)}{P(C \notin \D \mid E_1, I)}=\frac{P(S=C\mid E_1, I)}{P(S\neq C \mid E_1, I)}=\frac{1}{np}\frac{P(C \in \D \mid I)}{P( C \notin \D \mid I)},
\end{equation}
with corresponding likelihood ration $1/np$. If the database is comprised of individuals coming from different subpopulations, then \eqref{effodds} does not hold. However, in that case one may view the database as a disjoint union $\D=\D_1\cup \dots \cup D_m$, where $\D_m$ is the subset of $\D$ containing individuals from subpopulation $i$. For each of these separately, \eqref{effodds} holds.

It is rather interesting to see what happens with the odds on $S=C$ (given $E_1$ and $I$) when the size of the database grows.
It may seem from (\ref{effodds}) that as $n$ grows, the odds on $S=C$ decrease. However, this is not true in general, since $P(C \in \D \mid I)$ may also depend on $n$.
It does, however, mean that enlarging a database does not necessarily improve its effectiveness, in the sense of increasing the odds \eqref{effodds} on a unique match being with the true offender.
For example, suppose that a database $\D_n$ of size $n$ yields $P(C \in \D_n \mid I)=q_n$, and that a larger database $\D_{2n}$ of size $2n$ yields $P(C \in \D_{2n} \mid I)=q_{2n}$. If $\D_n \subset \D_{2n}$ then naturally $q_{2n} \geq q_n$, but the probability that $S=C$ given a unique match in $\D_{2n}$ is greater than the probability that $S=C$ given a unique match in $\D_n$ only when 
\[ \frac{q_{2n}}{1-q_{2n}}>2\frac{q_n}{1-q_n}.\]
This can be explained intuitively: if one adds many people who are unlikely to be $C$ to the database, then the probability of a chance match with one of these new individuals outweighs the fact that the probability that $C$ has been added to the database has increased in the sense that it becomes less likely that a unique match actually is a match with the criminal.

Hence the value of a unique match may increase or decrease with the size of the database, and it is not hard to see that the probability of a unique match itself may
(independently) decrease or increase. 

\subsection{Conclusions} 
\begin{itemize}
\item If it is known with whom the match is, say with $x_i$, then (cf. \eqref{k=1}) the posterior probability of guilt is given by
\[ \frac{P(C=x_i \mid E_\D, I)}{P(C \neq x_i \mid E_\D, I)}=\frac{\alpha_i}{p P(C \notin \D \mid I)}.\]
Notice that this quantity only depends on $\alpha_i=P(C=x_i\mid I)$, on the likelihood $p$ of a chance match with $x_i$ and on the a priori probability that the database contains the criminal. As the database increases, $P(C \notin \D \mid I)$ decreases but depending on $\alpha_i/p$ the posterior probability $P(C=x_i \mid E_\D, I)$ may be greater or smaller than for a smaller database.

\item
If it is not specified with which individual the match is, and the probability of having $\Gamma$ is $p$ for everyone in the database, then the posterior probability that the match is with the criminal is given by, cf. \eqref{effodds}, 
\[\frac{P(S=C\mid E_1, I)}{P(S\neq C \mid E_1, I)}=\frac{1}{np}\frac{P(C \in \D \mid I)}{P( C \notin \D \mid I)}.\]
These odds describe the long-term behaviour of the database, i.e., the proportion in the long run of unique matches that are matches with the true criminal. Naturally, enlarging the database always increases the probability that the criminal is contained in it. But the probability of a unique match may increase or decrease, and (independently) the value of a unique match may increase or decrease. In many cases, in an enlarged database the probability of a unique match increases, but the probability of a unique match being with the true offender decreases. 
\end{itemize}

\section{Examples}

In this section we illustrate the obtained results by considering some examples. We have chosen to cast most of these examples in a \textsc{dna}-setting, as this provides one of the few types of forensic evidence that are so well understood that more or less exact computations can be performed.

The uncertainty surrounding DNA-profile frequency estimates depends on the size of the database from which allele frequencies are estimated. A possible model is to define a prior distribution of allele frequencies, and to update this distribution with the database to obtain a posterior distribution. An often used approach is to use Dirichlet distributions (see \cite{TC} for an account of the method and a discussion on the sensitivity for the choice of prior). Doing this for a database containing alleles of 230 persons (for many forensic labs the actual size of their database is a few hundred individuals), it seems (based on simulations for DNA-profiles with six or seven loci and frequencies between $10^{-10}$ and $10^{-7}$) reasonable to use a standard deviation $p/3 \leq \sigma \leq 2p/3$ in the below examples.

We will in each example freely use the notation introduced in the section that it illustrates.
\subsection{Classical island problem with uncertain $\Gamma$-frequency}
We start with the simple version of a homogeneous population $X$ of size $N+1$ and profile frequency $p$. As we have seen (cf.\ \eqref{ptc} and \eqref{pte}), the posterior probability of guilt is equal to $P_s(G\mid I,E)=1/(1+N(p+\sigma^2/p))=1/(1+Np')$. With $p/3 \leq \sigma \leq p$, we get $p'\in [9p/8, 13p/9]$. Thus, the effect of the uncertainty about $p$ is to effectively increase $p$, or equivalently, to decrease the likelihood ratio associated to $I,E$ or to $E$. It may be prudent to use $\sigma=p$. For example, with $N=10^7, p=10^{-8}, \sigma=p$, we have $P_s(G\mid I,E)=0.83$ instead of $0.91$.

\subsection{Subpopulations and likelihood ratios}

We now illustrate the results of Section \ref{hetero}. Suppose that a crime has been committed in a heterogeneous population $X=X_1 \cup X_2$, with $N_1=10^7$ and $N_2=10^5$. Prior to DNA-analysis it is estimated that the crime could equally probably have been committed by a member of $X_1$ as by a member of $X_2$, i.e., $\beta_1=\beta_2=0.5$. Now a DNA-trace of the criminal is found, giving rise to a profile $\Gamma$. The forensic lab calculates $p_1=10^{-9}$ and $p_2=10^{-8}$. 

\subsubsection{Unconditioned on the profile} The likelihood ratio \eqref{LR3} (taking both the fact that the criminal and the suspect have $\Gamma$ as evidence) equals $1/(p_1\beta_1+p_2\beta_2)=1.8 \times 10^8$. This likelihood ratio holds for any suspect $s$, as long as $S$ is independent of $C$. 

With this likelihood ratio we obtain, for $s \in X_1$, posterior odds in favour of guilt equal to \[(p_1\beta_1+p_2\beta_2)^{-1}\beta_1/N_1\approx 9,\] corresponding to (cf.\ \eqref{solheteropop2}) $P_s(G\mid I,E)=0.9$. For $s \in X_2$ the posterior odds are \[(p_1\beta_1+p_2\beta_2)^{-1}\beta_2/N_2\approx 910,\] such that $P_s(G\mid I,E)=0.999$.

\subsubsection{Conditional on the profile} Given the fact that $C$ has $\Gamma$ and the frequencies $p_1, p_2$, we can also first calculate $P(C \in X_i \mid I)=\alpha_i$. This gives $\alpha_1=0.09$ and $\alpha_2=0.91$: since the profile $\Gamma$ is rarer in $X_1$, it is much more likely that the criminal is from $X_2$. The odds on $C$ belonging to $X_1$ are $\alpha_1/\alpha_2=10$. If this is taken as information relative to which everything else is conditioned, then the likelihood ratio associated to having $\Gamma$, is the inverse random match probability for the suspect: $1/p_1$ or $1/p_2$. This gives rise to the same $P_s(G\mid I,E)$: if $s \in X_1$ then the posterior odds are \[ p_1^{-1}\alpha_1/(N1-\alpha_1) \approx \alpha_1/(N_1p_1)=0.09/(10^{-9}10^7)=9,\] as above. Similarly, for $s \in X_2$, we get posterior odds \[ p_2^{-1}\alpha_2/(N2-\alpha_2) \approx \alpha_1/(N_2p_2)= 0.91/(10^{-8}10^5)=910,\] as above.

\subsubsection{Consequences of errors} When statements are made regarding the subpopulation to which $C$ belongs, one has to be careful to note whether or not $I$ has been taken into account. Indeed taking $\alpha_i$ equal to $\beta_i$, that is, $\alpha_1=\alpha_2=0.5$, we overestimate posterior odds in favour of guilt with a factor 10 for suspects from $X_1$ and underestimate them with the same factor for suspects from $X_2$. This is a serious overestimate of the actual odds for suspects from $X_1$. In this example, it leads to a posterior probability of guilt of $0.98$ (instead of $0.90$).

Finally, we note that if that the forensic lab assumes $p_2=p_1=10^{-9}$ for both populations, e.g. because it always uses the population frequencies of the dominant population $X_1$, then we arrive at $\alpha_i=\beta_i$. The posterior odds in favour of guilt will in that case be calculated to be $p_1^{-1}\alpha_i/(N_i-\alpha_i)\approx \alpha_i/(p_1N_i)$ for $s \in X_i$. In this example, these odds are 50 for $s \in X_1$ and 5000 for $s \in X_2$ which is an overestimate in both cases. 

\subsection{Subpopulations: general case} We next illustrate the results that we have obtained for the case where the populations is heterogeneous w.r.t.\ $\Gamma$-probability, and there is uncertainty about the profile frequency in each population, as well as uncertainty about the subpopulation to which an individual belongs. This is described in section \ref{uncertain}. Since there are many parameters that can be varied, we will keep some of them fixed throughout. We assume that the population consists of three disjoint subpopulations $X_1,X_2,X_3$, where $X_1$ is the dominant one, and the others are much smaller. We set $N_1=20\cdot 10^6, N_2=10^6, N_3=10^5$ and $\sigma=p/2$. We will compare the true posterior probability of guilt $P(G \mid I, E)$ with the probability obtained assuming that for $X_2, X_3$ the same $\Gamma$-frequency $p_1$ is used as for $X_1$. This allows one to judge what the consequences are of having a subpopulation without knowing so. For example, there may be a region of the country with a relatively high $\Gamma$-frequency due to its relative isolation in the past. In practice it can be difficult to say with certainty if a given individual belongs to that subpopulation.

We compute for several choices of $p_i$, $\epsilon_i$ and $\beta_i$ the true probability of guilt and compare it to what one would obtain if $p_2,p_3$ would be ignored, namely \eqref{ptc} with $N+1=N_1+N_2+N_3$ and $p=p_1$. We denote this result with $P^{{\rm hom}}(G\mid I,E)$ and call it the {\it naive} probability of guilt. 

\begin{example}\label{exprobs1} {\rm Let $p_1=10^{-8}, p_2=10^{-7}, p_3=10^{-6}$. We keep the $\epsilon_i$ fixed to a choice where it is 90\% certain that $S \in X_1$, not knowing $I$ or $E$. The results are summarized in Table \ref{probs1}.
\begin{table}[h]\caption{Guilt probabilities for $p_1=10^{-8}, p_2=10^{-7}, p_3=10^{-6}$}\label{probs1}\begin{tabular}{|c|c|c|c|}\hline  $(\epsilon_1,\epsilon_2,\epsilon_3)$ & $(\beta_1,\beta_2,\beta_3)$ & $P(G|I,E)$ & $P^{{\rm hom}}(G|I,E)$ \\ \hline  (0.9,0.05,0.05) & (0.999,0.0005,0.0005) & 0.50 & 0.79  \\  (0.9,0.05,0.05) & uniform & 0.70 & 0.79 \\ (0.9,0.05,0.05) & (0.99,0.005,0.005) & 0.74 & 0.79  \\  (0.9,0.05,0.05) & (0.9,0.05,0.05) & 0.84 & 0.79 \\ \hline \end{tabular} \end{table}
Notice that the true probability of guilt may be smaller or greater than the naive probability. In the first line with $\beta_1=0.999$, there is considerable uncertainty as to the subpopulation to which $S$ belongs given $I, E$; in fact $P_s(S \in X_1 \mid I, E)=0.40, P_s(S \in X_3 \mid I,E)=0.56$. Since for this choice of parameters $P_s(G \mid I, E, S \in X_3)$ (given by \eqref{post}) is only  0.32, we get a probability of guilt equal to 0.50, much smaller than the naive probability. However, as $\beta_1$ decreases, so does $P_s(S \in X_1 \mid I, E)$, and $P_s(S \in X_3 \mid I,E)$ grows. In the last line of Table \ref{probs1}, $P_s(S \in X_3 \mid I, E)$ is large (equal to 0.95), so the posterior probability of guilt is predominantly given by \eqref{post} applied to $s \in X_3$, which is 0.89 for these parameters.}
\end{example}

\begin{example}\label{exprobs2} {\rm As observed above, we obtain the same probabilities $P(G \mid I,E)$ (and of course, the same naive probability of guilt), when in the above example $\epsilon$ and $\beta$ are exchanged.
The explanation for these probabilities is somewhat different. In the first line of Table \ref{probs1} (now with $\epsilon_1=0.999$), it is quite likely that $S$ belongs to $X_1$ given $I, E$; in fact $P_s(S \in X_1 \mid I, E)=0.79, P_s(S \in X_3 \mid I,E)=0.21$. Since for this choice of parameters $P_s(G \mid I, E, S \in X_1)$ (given by \eqref{post}) is only  0.40, we get a probability of guilt equal to 0.50, much smaller than the naive probability. However, exactly as for Example \ref{exprobs1}, as $\epsilon_1$ decreases, so does $P_s(S \in X_1 \mid I, E)$, and $P_s(S \in X_3 \mid I,E)$ grows.} 
\end{example}

These examples show that the effect of having subpopulations can be considerable when the profile is more common among the smaller subpopulations, even when both $S$ and $C$ are likely issued from the largest subpopulation. The magnitude and the direction of the subpopulation effect depend strongly on the a priori probabilities for $S$ and $C$ to belong to each of the subpopulations.

\begin{example} {\rm Letting $S$ and $C$ be likely issued from $X_2$ or $X_3$, we get a posterior probability of guilt between 0.80 and 0.85 which does not depend strongly on the precise choice of $\epsilon_i$ and $\beta_i$. This is understandable since these choices all make $P(S \in X_1 \mid I,E)$ small, and \eqref{post} applied to $X_2$ and $X_3$ yields 0.89 for both populations (note that they have the same expected number of $\Gamma$-bearers).}
\end{example}
\begin{example}
{\rm Consider the case where $p_2$ and $p_3$ are smaller than $p_1$, for example $p_1=10^{-8}, p_2=10^{-9}, p_3=10^{-10}$. The population as a whole then has a smaller number of expected $\Gamma$-bearers compared to when $p_2=p_3=p_1$. The true probability of guilt exceeds the naive probability unless one is almost sure that $S$ and $C$ are from different subpopulations, as illustrated in Table \ref{probs3}.
\begin{table}[h]\caption{Guilt probabilities for $p_1=10^{-8}, p_2=10^{-9}, p_3=10^{-10}$}\label{probs3}\begin{tabular}{|c|c|c|c|}\hline  $(\epsilon_1,\epsilon_2,\epsilon_3)$ & $(\beta_1,\beta_2,\beta_3)$ & $P(G|I,E)$ & $P^{{\rm hom}}(G|I,E)$ \\ \hline  uniform & uniform & 0.80 & 0.79 \\ (0.9,0.05,0.01) & (0.9,0.05,0.05) & 0.80 & 0.79  \\  (0.2,0.6,0.2) & (0.2,0.6,0.2) & 0.98 & 0.79 \\ (0.1,0.3,0.6) & (0.3,0.3,0.4) & 0.97 & 0.79  \\  (0.9,0.09,0.01) & (0.01,0.01,0.98) & 0.88 & 0.79 \\ (0.99,0.009,0.001) & (0.001,0.001,0.998) & 0.57 & 0.79\\ \hline \end{tabular} \end{table}}
\end{example}

\subsection{$\Gamma$-correlation: relatedness}\label{excorr} Suppose that $C$ has DNA-profile $\Gamma$ with a population frequency of $10^{-7}$, i.e., $p_x=10^{-7}$ for all $x \in X$. Now we select $s$ from $X$, and $\Gamma_s=1$. Suppose that $X=\{s,y_1,y_2,y_3,z_1,\dots,z_N\}$ and $N=10^6$, such that $c_{y_i,s}=10^{-3}$ and $c_{z_i,s}=p_{z_i}=10^{-7}$. Here we model a situation in which the suspect $s$ has three brothers, whose $\Gamma$-probability is $10^{-3}$ given $\Gamma_s=1$, and that the rest of the population is unrelated to $s$. If the a priori probabilities are $P_s(C=s\mid I)=0.4, P_s(C=y_i\mid I)=0.1, P_s(C=z_j\mid I)=0.3/10^6$ then the correction factor \eqref{corrfactor} for the $\Gamma$-correlation is equal to
\[ \frac{0.6}{3\cdot 10^4 \cdot 0.1 + 10^6\cdot 1 \cdot 0.3\cdot 10^{-6}}=\frac{0.6}{0.3+3000}\approx \frac{1}{5000},\]
meaning that the likelihood ratio associated to $\Gamma_s=1$ has been made 5000 times smaller, reducing it from $1/p_s=10^7$ to 2000.

For the posterior probability of guilt $P_s(G \mid I, E)$, this means that it is reduced from 
\[ \frac{\frac{2}{3}10^7}{1+\frac{2}{3}10^7} \approx 1-\frac{3}{2}10^{-7}\]
that we would obtain without $\Gamma$-correlation, to approximately
 1-3/4000.
\subsection{Biased search} \label{exsel} We recast the example given in Section \ref{excorr} in the setting of a biased search, to demonstrate the equivalence noted in Section \ref{equiv}. As in \ref{excorr}, $p_x=10^{-7}$ for all $x \in X$, we suppose that $S=s$ has been selected and that $\Gamma_s=1$,and that there are $y_1,y_2,y_3 \in X$ such that $\sigma_{x,y_i}=10^4\sigma_{x,x}$ for $i=1,2,3$ and $\sigma_{x,z_i}=\sigma_{x,x}$ for all $i$. The prior odds are as in 
Section \ref{excorr}. Then the likelihood ratio associated to the evidence $\Gamma_s=1$ is reduced by a factor of about 5000, as in Example \ref{excorr}. In that example, the value was decreased since finding $\Gamma$ in $s$ made it more probable that population members $y_1,y_2,y_3$, which have non-negligible prior probabilities of guilt, also have $\Gamma$. In this situation, it is due to the fact that the selection procedure is such that if $s$ is selected, it becomes less likely that $s$ is guilty:
\[ P(C=s \mid S=s, I)=\frac{P(C=s \mid I)}{\sum_{y \in X}\frac{\sigma_{s,y}}{\sigma_{s,s}}P(C=y \mid I)} 
\approx \frac{1}{7500},\]
which is considerably less than $0.4$.
The fact that $s$ has $\Gamma$ then raises the probability of guilt to
approximately $1-3/4000$ as above.

\subsection{Database effectiveness} In \eqref{effodds} we have computed the odds in favour of a unique database match being with the true criminal. If the database is a random sample of the population in the sense that $P(C \in \D \mid I)=|\D|/|X|=n/N$, then this equation reads
\[ \frac{P(S=C \mid E_1, I)}{P(S\neq C\mid E_1, I)}=\frac{1}{p(N-n)},\] which is monotonically increasing in $n$, going from $1/((N-1)p)$ for $n=1$ to infinity for $n=N$. It is not hard to derive this directly: since $n-1$ persons have been shown not to possess $\Gamma$, the population that can not be excluded has size $N-n+1$. In that population, only the $\Gamma$-status of one individual (the one that matched in $\D$) is known. Since $\D$ was a random sample as defined above, the classical solution \eqref{sol} applies.

If the database is not a random sample from the population in the above sense, then the situation is more interesting and quite different. 

\begin{example} {\rm Let $p=10^{-7}$ and suppose that with $n=10^5$ one has $P(C\in \D \mid I)=0.2$. For example, this may be because the database consists of previously convicted individuals and based on the probability of a rightful conviction and of recidivism one arrives at such an estimate. For database $\D$, the odds that a unique match is with $C$ are 25 to one, or equivalently, $P(S=C \mid I, E_1)=0.96$. 

It may be possible to enlarge $\D$ to $\D'$ with $|\D'|=n'$ such that $P(S=C \mid I,E_1)=0.5$, but only at the cost of adding very many individuals into $\D'$, e.g. with $n'=2\cdot 10^6$. In that case, the odds \eqref{effodds} on a unique match being with the offender decrease to 5, i.e., one in six of such matches will be with an innocent person.

The probability of actually obtaining a unique match is given by
\[ P(E_1)=P(C \in \D \mid I)(1-p)^n+P(C \notin \D \mid I)np(1-p)^{n-1}.\]\
For database $\D$, this evaluates to 0.206 and for database $\D'$ to 0.491\%. Thus, in $\D'$ a search with a DNA-profile with population frequency $10^{-7}$ will yield a unique match about half of the time, but only 5 out of 6 of these will be with the true offender. About 10\% of such searches will result in two or more hits, and about 40\% will not result in any hit.

When multiple matches are found, it is more likely that one of them is with the true offender but not a near certainty: e.g., in case two matches are found (which happens with probability 0.09), about one in ten of such double matches are both coincidental.

For the original database $\D$, about 20,6\% of searches result in a unique match, almost all of which are with the offender; in the remaining cases one almost always has no hits: the probability of having more than one match being 0.002.}
\end{example}

\begin{example} {\rm Suppose that the database is set up and expanded such that if it has size $n$ then $P(C \in \D)=\sqrt{n/N}$. This is a model for a database in which individuals with higher prior probability of guilt are put in the database with higher probability. For example, if $\D$ contains the DNA-profiles of 10\% of the population, then it contains $C$ with probability 0.31. If $\D$ is enlarged to contain 30\% of the population, then it contains $C$ with probability 0.55.

In that case, the odds on a unique match being with the criminal in a database of size $n$ are minimal for $n=N/4$. An example for $N=2\cdot10^7, p=10^{-8}$ is given in Figure \ref{odds}. With $n=N/4$ the odds in favour of a unique match being with the true offender are 20. As the plots show, when the database is relatively small the odds on a match being with the true offender decrease rapidly, e.g. from 105 if $n=50.000$ to 55 if $n=200.000$. As $n$ grows further, the odds decrease (slowly) to 20 for $n=5\cdot 10^6=N/4$. When $n$ grows further, the odds increase again. When 50\% of the population is included ($n=10^7$), they are 24.
\begin{figure}[h]
\caption{Database effectiveness with $p=10^{-8}, N=2\cdot 10^7, P(C \in \D)=\sqrt{n/N}$}
\label{odds}
\centerline{\includegraphics[width=7cm]{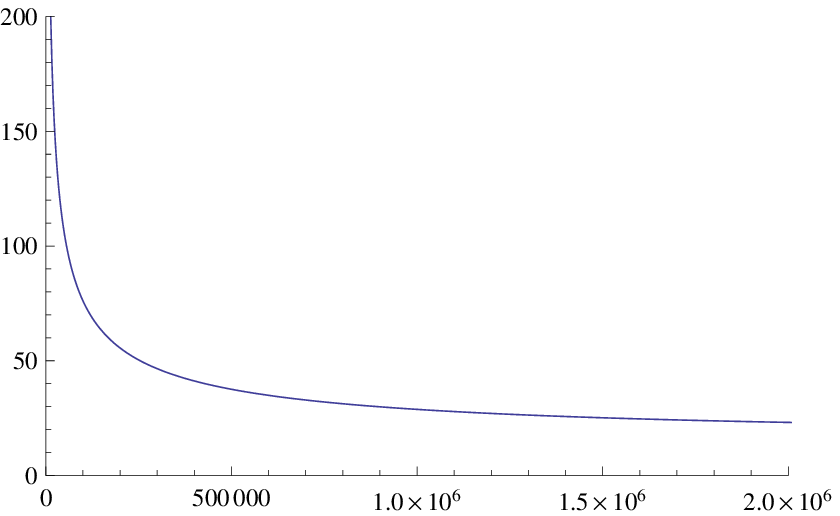}\includegraphics[width=7cm]{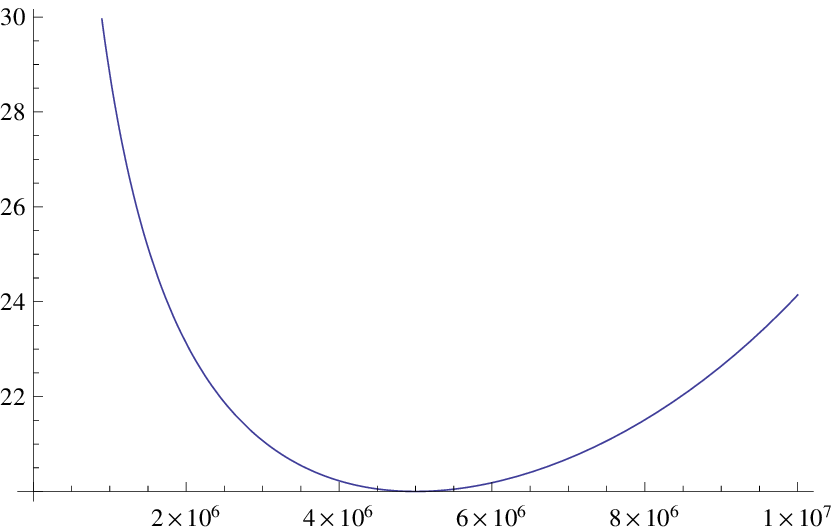}}
\end{figure}
}
\end{example}

Thus, enlarging a database may at the same time increase the chance of obtaining a unique match from it, and diminish the value of such a match in the sense that the probability of it being with the true offender decreases. These examples suggest that the idea that the larger the database, the better, needs to be put into perspective. It is of course true that enlarging a database increases the probability that the criminal is included. It is also obvious that given a unique match in the database, the probability that it is with the criminal increases when the database is expanded and does not yield additional matches. But as we have seen, it does not follow that hits in larger databases are stronger evidence for guilt than hits in smaller databases.

\end{document}